
\input amssym  

\hsize=15truecm
\hoffset=.46truecm
\vsize=23.7truecm
\voffset=.46truecm

\ifx\pdfpagewidth\undefined

\else
\pdfpagewidth=210truemm
\pdfpageheight=297truemm
\fi

\newdimen\normalparindent

\iffalse  

\normalparindent=24pt

\font\elevenrm=cmr10 at 11pt 
\font\eightrm=cmr8 
\font\sixrm=cmr6 

\font\eleveni=cmmi10 at 11pt 
\font\eighti=cmmi8
\font\sixi=cmmi6

\font\elevensy=cmsy10 at 11pt 
\font\eightsy=cmsy8
\font\sixsy=cmsy6

\font\elevenex=cmex10 at 11pt 

\font\elevenbf=cmbx10 at 11pt 
\font\eightbf=cmbx8
\font\sixbf=cmbx6

\font\eleventt=cmtt10 at 11pt 
\font\elevensl=cmsl10 at 11pt 
\font\elevenit=cmti10 at 11pt 

\textfont0=\elevenrm \scriptfont0=\eightrm \scriptscriptfont0=\sixrm
\def\rm{\fam0\elevenrm}
\textfont1=\eleveni \scriptfont1=\eighti \scriptscriptfont1=\sixi
 
\textfont2=\elevensy \scriptfont2=\eightsy \scriptscriptfont2=\sixsy
\def\cal{\fam2}
\textfont3=\elevenex \scriptfont3=\elevenex \scriptscriptfont3=\elevenex
\textfont\itfam=\elevenit
\def\it{\fam\itfam\elevenit}
\textfont\slfam=\elevensl
\def\sl{\fam\slfam\elevensl}
\textfont\bffam=\elevenbf \scriptfont\bffam=\eightbf
\scriptscriptfont\bffam=\sixbf
\def\bf{\fam\bffam\elevenbf}
\textfont\ttfam=\eleventt
\def\tt{\fam\ttfam\eleventt}

\skewchar\eleveni='177 \skewchar\eighti='177 \skewchar\sixi='177
\skewchar\elevensy='60 \skewchar\eightsy='60 \skewchar\sixsy='60

\font\sc=cmcsc10 at 11pt 

\font\cmcyr=cmcyr10 at 11pt
\font\cmcti=cmcti10 at 11pt
\font\cmccsc=cmccsc10 at 11pt

\smallskipamount=3.5pt plus 1pt minus 1pt
\medskipamount=7pt plus 2pt minus 2pt
\bigskipamount=14pt plus 2pt minus 2pt
\normalbaselineskip=14pt
\normallineskip=1pt
\normallineskiplimit=0pt
\jot=3.5pt

\normalbaselines
\rm

\else  

\normalparindent=20pt

\font\sc=cmcsc10

\font\cmcyr=cmcyr10
\font\cmcti=cmcti10
\font\cmccsc=cmccsc10

\fi

\def\cyrchardefs{%
\chardef\yo=60
\chardef\Yo=62
\chardef\yu=64
\chardef\a=65
\chardef\b=66
\chardef\ts=67
\chardef\d=68
\chardef\ye=69
\chardef\f=70
\chardef\g=71
\chardef\kh=72
\chardef\i=73
\chardef\j=74
\chardef\k=75
\chardef\l=76
\chardef\m=77
\chardef\n=78
\chardef\o=79
\chardef\p=80
\chardef\ya=81
\chardef\r=82
\chardef\s=83
\chardef\t=84
\chardef\u=85
\chardef\zh=86
\chardef\v=87
\chardef\soft=88
\chardef\y=89
\chardef\z=90
\chardef\sh=91
\chardef\e=92
\chardef\shch=93
\chardef\ch=94
\chardef\hard=95
\chardef\Yu=96
\chardef\A=97
\chardef\B=98
\chardef\Ts=99
\chardef\D=100
\chardef\Ye=101
\chardef\F=102
\chardef\G=103
\chardef\H=104
\chardef\I=105
\chardef\J=106
\chardef\K=107
\chardef\L=108
\chardef\M=109
\chardef\N=110
\chardef\O=111
\chardef\P=112
\chardef\Ya=113
\chardef\R=114
\chardef\S=115
\chardef\T=116
\chardef\V=119
\chardef\Soft=120
\chardef\Y=121
\chardef\Z=122
\chardef\Sh=123
\chardef\E=124
\chardef\Shch=125
\chardef\Ch=126
\chardef\Hard=127
}
\def\cyr{\cmcyr\cyrchardefs}
\def\cycsc{\cmccsc\cyrchardefs}
\def\cyti{\cmcti\cyrchardefs}

\def\S{\mathhexbox278\thinspace}

\def\square{\hbox to.77778em{%
\hfil\vrule\vbox to.675em{\hrule width.6em\vfil\hrule}\vrule\hfil}}

\long\def\remerciements#1\par{\medbreak\noindent{\it
    Remerciements\/}.\enspace #1\par\medbreak}
\def\definition#1\par{\medbreak\noindent{\bf D\'efinition.}\enspace
  #1\par\medbreak}
\def\exemple#1\par{\medbreak\noindent{\bf Exemple.}\enspace
  #1\par\medbreak}
\long\def\remarque#1\par{\medbreak\noindent{\it Remarque\/}.\enspace
#1\par\medbreak}
\long\def\remarques#1\par{\medbreak\noindent{\it Remarques\/}.\enspace
#1\par\medbreak}
\def\exercice#1\par{\medbreak\noindent{\bf Exercice.}\enspace
#1\par\medbreak}
\def\notation#1\par{\medbreak\noindent{\bf Notation.}\enspace
#1\par\medbreak}
\def\preuve{\noindent{\it Preuve\/}.\enspace}
\def\preuvefinie{\nobreak\hfill\quad\square\par\medbreak}

\def\lineover#1{{\offinterlineskip\mathchoice
{\setbox0=\hbox{$\displaystyle#1$}%
\vbox{\kern .1ex\hbox to\wd0{\kern .2ex\leaders\hrule height .1ex%
\hfill\kern .2ex}\kern .1ex\box0}}
{\setbox0=\hbox{$\textstyle#1$}%
\vbox{\kern .1ex\hbox to\wd0{\kern .2ex\leaders\hrule height .1ex%
\hfill\kern .2ex}\kern .1ex\box0}}
{\setbox0=\hbox{$\scriptstyle#1$}%
\vbox{\kern .08ex\hbox to\wd0{\kern .16ex\leaders\hrule height .08ex%
\hfill\kern .16ex}\kern .08ex\box0}}
{\setbox0=\hbox{$\scriptscriptstyle#1$}%
\vbox{\kern .06ex\hbox to\wd0{\kern .12ex\leaders\hrule height .06ex%
\hfill\kern .12ex}\kern .06ex\box0}}}}

\def\isomorphism#1{\mathrel{\mathop{\longrightarrow}%
\limits^{#1}_{\raise0.5ex\hbox{$\scriptstyle\sim$}}}}

\def\injlim{\mathop{\vtop{\offinterlineskip\halign{##\cr
 \hfil\rm lim\hfil\cr\noalign{\kern.1ex}\rightarrowfill\cr
 \noalign{\kern-.4ex}\cr}}}}
\def\projlim{\mathop{\vtop{\offinterlineskip\halign{##\cr
 \hfil\rm lim\hfil\cr\noalign{\kern.1ex}\leftarrowfill\cr
 \noalign{\kern-.4ex}\cr}}}}

\def\blank{\mkern12mu}

\def\textfrac#1/#2{{\textstyle{#1\over#2}}}

\def\relativediag#1#2#3#4#5#6{\vcenter{\baselineskip=3ex \halign{
\hfil$##$&$##$&$##$\hfil\cr
#1\quad& \hfilneg\buildrel#2\over\longrightarrow\hfilneg& \quad#3\cr
\lower1ex\llap{$\scriptstyle#4\hskip-1ex$}\searrow& \quad&
\swarrow\lower1ex\rlap{$\hskip-1ex\scriptstyle#5$} \cr
& \hfilneg#6\hfilneg&\cr}}}

\def\trianglediag#1#2#3#4#5#6{\vcenter{\baselineskip=3ex \halign{
\hfil$##$&$##$&$##$\hfil\cr
#1\quad& \hfilneg\buildrel#2\over\longrightarrow\hfilneg& \quad#3\cr
\lower1ex\llap{$\scriptstyle#6\hskip-1ex$}\nwarrow& \quad&
\swarrow\lower1ex\rlap{$\hskip-1ex\scriptstyle#4$} \cr
& \hfilneg#5\hfilneg&\cr}}}

\def\correspondence#1#2#3#4#5{\vcenter{\baselineskip=3ex \halign{
\hfil$##$&$##$&$##$\hfil\cr
&\hfilneg#1\hfilneg\cr
\raise1ex\llap{$\scriptstyle#2$}\swarrow&&\searrow
\raise1ex\rlap{$\scriptstyle#3$}\cr
#4&&#5\cr}}}

\newif\iffirstpar
\everypar{\iffirstpar\parindent=\normalparindent\firstparfalse\fi}

\def\sectionheading#1{\subcount=0 \subsectioncount=0 \eqcount=0
  \bigskip\vskip\parskip
  \leftline{\bf #1}\nobreak\smallskip\firstpartrue\parindent=0pt}

\newif\ifappendix \appendixfalse

\def\currentsection{\ifappendix A\else\number\sectioncount\fi}

\def\section#1\par{\advance\sectioncount by1%
  \edef\currentlabel{\currentsection}%
  \sectionheading{\currentsection.\enspace#1}}

\def\unnumberedsection#1\par{\sectionheading{#1}}

\def\references#1{
\unnumberedsection#1\par
\normalparindent=25pt
\parindent=\normalparindent
\parskip=1ex plus 0.5ex minus 0.2ex}

\def\subsection#1\par{\medbreak\penalty-200\advance\subsectioncount by1%
  \edef\currentlabel{\currentsection.\number\subsectioncount}%
  \leftline{\it\currentsection.\number\subsectioncount.\enspace#1}%
  \smallskip\parindent=0pt\firstpartrue}

\newwrite\auxfile

\newcount\sectioncount \sectioncount=0
\newcount\subsectioncount 
\newcount\subcount 
\newcount\eqcount 

\def\subno{\global\advance\subcount by1\relax
  \currentsection.\number\subcount
  \xdef\currentlabel{\currentsection.\number\subcount}}
\def\proclaim #1. #2\par{\medbreak
  \noindent{\bf#1~\subno.\enspace}{\sl#2\par}%
  \ifdim\lastskip<\medskipamount \removelastskip\penalty55\medskip\fi}
\def\proclaimx #1 (#2). #3\par{\medbreak
  \noindent{\bf#1~\subno\ \rm (#2).\enspace}{\sl#3\par}%
  \ifdim\lastskip<\medskipamount \removelastskip\penalty55\medskip\fi}

\newdimen\algindent
\def\plusindent{\advance\algindent by \parindent}
\def\minusindent{\advance\algindent by-\parindent}


\newcount\algstepcount

\long\def\algorithm (#1). #2\endalgorithm{\medbreak
  \algindent=0pt%
  \algstepcount=0%
  \noindent{\bf Algorithme~\subno} (#1). {\sl#2}\par\medbreak}

\def\step{\advance\algstepcount by1
\edef\currentlabel{\number\algstepcount}
\smallskip\hangindent\parindent
\advance\hangindent by\algindent\indent
\llap{{\bf \the\algstepcount.}\enspace}\kern\algindent
\ignorespaces}


\def\labeldef#1#2{\expandafter\gdef\csname L@#1\endcsname{#2}}
\def\label#1{%
  \expandafter\xdef\csname L@#1\endcsname{\currentlabel}%
  \write\auxfile{\string\labeldef{#1}{\csname L@#1\endcsname}}%
  \ignorespaces}
\def\ref#1{{\rm\expandafter\ifx\csname L@#1\endcsname\relax
  \message{Undefined label `#1'}??\else
  \csname L@#1\endcsname\fi}}

\def\eqdef#1#2{\expandafter\gdef\csname E@#1\endcsname{#2}}
\def\eqnumber#1{\global\advance\eqcount by1\relax
  \eqno{\rm(\currentsection.\number\eqcount)}%
  \expandafter\xdef\csname E@#1\endcsname{%
    \currentsection.\number\eqcount}%
  \write\auxfile{\string\eqdef{#1}{\csname E@#1\endcsname}}}
\def\eqref#1{{\rm(\expandafter\ifx\csname E@#1\endcsname\relax
  \message{Undefined equation `#1'}??\else
  \csname E@#1\endcsname\fi)}}

\newcount\refcount \refcount=0
\def\citedef#1#2{\expandafter\gdef\csname C@#1\endcsname{#2}}
\def\cite#1{\expandafter\ifx\csname C@#1\endcsname\relax
  \message{Undefined reference `#1'}\citedef{#1}{??}\fi
  \expandafter\gdef\csname R@#1\endcsname{\relax}%
  [\csname C@#1\endcsname]}
\def\citex#1#2{\expandafter\ifx\csname C@#1\endcsname\relax
  \message{Undefined reference `#1'}\citedef{#1}{??}\fi
  \expandafter\gdef\csname R@#1\endcsname{\relax}%
  [\csname C@#1\endcsname, #2]}
\def\reference#1{\advance\refcount by 1%
  \expandafter\ifx\csname R@#1\endcsname\relax
  \message{Warning: reference `#1' not used}\fi
  \expandafter\edef\csname C@#1\endcsname{\the\refcount}%
  \write\auxfile{\string\citedef{#1}{\csname C@#1\endcsname}}%
  \item{[\csname C@#1\endcsname]}}

\newif\ifauxexists
\immediate\openin0=\jobname.aux
\ifeof 0
  \auxexistsfalse
\else
  \auxexiststrue
\fi
\immediate\closein0
\ifauxexists
  \input \jobname.aux
\else
  \message{No file `\jobname.aux'}
\fi
\openout\auxfile=\jobname.aux

\catcode`@=11
\def\f@encoding{OTF1}
\input francais.sty
\catcode`@=12

\selectlanguage{french}

\font\titlefont=cmssbx10 scaled \magstep2
\font\tenbfscaled=cmbx10 scaled \magstep2
\font\sevenbfscaled=cmbx7 scaled \magstep2
\font\fivebfscaled=cmbx5 scaled \magstep2

\newfam \titlefam
\textfont\titlefam=\tenbfscaled
\scriptfont\titlefam=\sevenbfscaled
\scriptscriptfont\titlefam=\fivebfscaled

\def\C{{\bf C}}
\def\mucan_#1{\mu^{\rm can}_{#1}}
\def\Nm{\mathop{\rm Nm}\nolimits}  
\def\Norm{{\rm N}}
\def\O{{\cal O}}
\def\P{{\bf P}}
\def\fp{{\frak p}}
\def\Pl{\Omega}  
\def\Q{{\bf Q}}
\def\Qbar{\lineover{\bf Q}}
\def\R{{\bf R}}
\def\vol{\mathop{\rm vol}\nolimits}
\def\Z{{\bf Z}}

\mathchardef\idot="213A
\def\fin{{\rm fin}}
\def\infin{{\rm inf}}
\def\stable{{\rm stable}}


{\titlefont\def\Qbar{\lineover{\fam\titlefam Q}}%
\centerline{Bornes optimales pour la diff\'erence entre la
  hauteur de Weil}
\smallskip
\centerline{et la hauteur de N\'eron--Tate sur les courbes
  elliptiques sur $\Qbar$}}
\medskip
\centerline{Peter Bruin}
\bigskip

\vfootnote{}{{\it Mathematics Subject Classification\/} (2010):
11G05,  
11G50,  
11Y35   

L'auteur a b\'en\'efici\'e du soutien financier du Fonds national
suisse (subsides nos.\ 124737 et 137920) et de l'hospitalit\'e du
Max-Planck-Institut f\"ur Mathematik, Bonn.}

\bigskip

{\narrower\narrower
\noindent{\it R\'esum\'e\/.}\enspace
Nous donnons un algorithme qui, \'etant donn\'ee une courbe elliptique
$E$ sur~$\Qbar$ sous la forme de Weierstra\ss, calcule l'infimum et le
supremum de la diff\'erence entre la hauteur na{\"\i}ve et la hauteur
canonique sur~$E(\Qbar)$.

\medskip
\noindent{\it Abstract\/.}\enspace 
We give an algorithm that, given an elliptic curve $E$ over~$\Qbar$ in
Weierstra\ss\ form, computes the infimum and supremum of the
difference between the na{\"\i}ve and canonical height functions
on~$E(\Qbar)$.\par}

\section Introduction

Soit $E$ une courbe elliptique sur~$\Qbar$ donn\'ee par une \'equation
de Weierstra\ss:
$$
E:y^2+a_1xy+a_3y=x^3+a_2x^2+a_4x+a_6.
\eqnumber{weierstrass}
$$
On supposera toujours que {\it les $a_i$ soient des entiers
alg\'ebriques\/}.

On s'int\'eresse \`a la hauteur de Weil (ou hauteur na{\"\i}ve) et
la hauteur de N\'eron--Tate (ou hauteur canonique)
$$
h,\hat h:E(\Qbar)\to\R
$$
(voir ci-dessous pour les normalisations).  La diff\'erence $h-\hat h$
est born\'ee sur~$E(\Qbar)$.  Il est donc naturel de se demander
si, \'etant donn\'ee une \'equation \eqref{weierstrass}, on peut
calculer le supremum et l'infimum de cette diff\'erence.

Des r\'esultats dans cette direction ont \'et\'e obtenus par
Dem$'$janenko \cite{Demjanenko} et Zimmer \cite{Zimmer},
Silverman \cite{Silverman}, Siksek \cite{Siksek}, Cremona, Prickett et
Siksek \cite{Cremona-Prickett-Siksek}, et Uchida \cite{Uchida}.  Le
lecteur est renvoy\'e \`a l'introduction
de \cite{Cremona-Prickett-Siksek} pour plus de d\'etails.

\subsection Notations

Dans cet article, \'etant donn\'e un corps de nombres $K$, on fixe les
notations suivantes:
\smallskip
\settabs\+\quad&$\Pl_K$, $\Pl^\fin_K$, $\Pl^\infin_K$\quad& bla\cr
\+& $\Pl_K$, $\Pl^\fin_K$, $\Pl^\infin_K$& l'ensemble des places
de~$K$, resp.\ des places finies, resp.\ des places infinies\cr
\+& $\Z_K$& l'anneau des entiers de~$K$\cr
\smallskip\noindent
Pour chaque place $v$:
\smallskip
\+& $|\blank|_v$& la valeur absolue normalis\'ee sur~$K$
correspondant \`a~$v$\cr
\+& $K_v$& le compl\'et\'e $v$-adique de~$K$\cr
\+& $\epsilon_v$& le degr\'e local $[K_v:\Q_p]$, o\`u $p$ est la place
de~$\Q$ au-dessous de~$v$\cr
\smallskip\noindent
Pour chaque place finie $v$:
\smallskip
\+& $\fp_v$& l'id\'eal maximal $\{x\in\Z_K\mid |x|_v<1\}$ de~$\Z_K$\cr
\+& $k_v$& le corps r\'esiduel $\Z_K/\fp_v$\cr
\smallskip\noindent
Soit $E$ est une courbe elliptique sur~$\Qbar$ donn\'ee par
une \'equation~\eqref{weierstrass}.  On d\'efinit les coefficients
$b_2$, $b_4$, $b_6$, $c_4$, $c_6$, le discriminant $\Delta_E$ et
l'invariant $j_E$ par les formules usuelles; voir par example
Tate \citex{Tate}{\S2}.

\subsection Hauteurs

La hauteur (logarithmique normalis\'ee) d'un point
$x=(x_0:\ldots:x_n)\in\P^n(\Qbar)$ est d\'efinie comme suit: soit
$K\subset\Qbar$ un corps de nombres contenant les $x_i$, alors
$$
h_{\P^n}(x)={1\over[K:\Q]}\sum_{v\in\Pl_K}
\log\max\{|x_0|_v,\ldots,|x_n|_v\}.
$$
On sait que le membre de droite est invariant par extension du
corps~$K$, de sorte que $h_{\P^n}$ est une fonction bien d\'efinie
sur~$\P^n(\Qbar)$.

Soit $E$ une courbe elliptique sur~$\Qbar$ donn\'ee par une \'equation
de Weierstra\ss~\eqref{weierstrass}, et soit $P\in E(\Qbar)$.  La {\it
hauteur de Weil\/} (ou {\it hauteur na{\"\i}ve\/}) de~$P$ est
d\'efinie par
$$
h(P)=h_{\P^1}(x(P)).
$$
La {\it hauteur de N\'eron--Tate\/} (ou {\it hauteur canonique\/})
de~$P$ est d\'efinie par
$$
\hat h(P)=\lim_{n\to\infty}n^{-2}h(nP).
$$
Notre fonction $\hat h$ co{\"\i}ncide avec celle utilis\'ee dans
l'article de Cremona, Prickett et
Siksek \cite{Cremona-Prickett-Siksek}; elle est double de celle
utilis\'ee dans l'article \cite{Silverman} et le
livre~\cite{Silverman, ATAEC} de Silverman.

\section La diff\'erence $h-\hat h$

Soit $E$ une courbe elliptique sur~$\Qbar$ donn\'ee par une \'equation
de Weierstra\ss~\eqref{weierstrass}, et soit $P\in E(\Qbar)$.

\subsection Hauteurs locales

Soit $K\subset\Qbar$ un corps de nombres tel que $E$ soit d\'efinie
sur~$K$.  Alors on a des hauteurs locales
$$
\lambda_v : E(K_v)\to\R\quad\hbox{pour tout }v\in\Pl_K
$$
et une d\'ecomposition de $\hat h$ en termes locaux
$$
[K:\Q]\hat h(P)=2\sum_{v\in\Pl_K}\epsilon_v\lambda_v(P)
\quad\hbox{pour tout }P\in E(K).
$$
La normalisation des $\lambda_v$ que nous utiliserons est celle du
livre de Silverman \citex{Silverman, ATAEC}{Chapter~VI}.

Pour $v\in\Pl_K$ et $P\in E(K_v)$, notons
$$
\phi_v(P)=\epsilon_v^{-1}\log\max\bigl\{1,|x(P)|_v\bigr\}-2\lambda_v(P).
$$
Alors on a
$$
\eqalign{
[K:\Q]\bigl(h(P)-\hat h(P)\bigr)&=\sum_{v\in\Pl_K}
\left(\log\max\bigl\{1,|x(P)|_v\bigr\}
-2\epsilon_v\lambda_v(P)\right)\cr
&=\sum_{v\in\Pl_K}\epsilon_v\phi_v(P).}
$$

\subsection Le discriminant stable

Quitte \`a \'elargir $K$, on peut supposer que $E$ ait r\'eduction
semi-stable.  Pour toute place finie $v$ de~$K$, on note $n_v$ le
nombre de composantes g\'eom\'etriques irr\'eductibles de la
r\'eduction de~$E$ modulo~$v$; cette r\'eduction est donc un
$n_v$-gone.  On note $\Delta_{E/K}^{\min}$ le discriminant minimal
de~$E$ sur~$K$, c'est-\`a-dire l'id\'eal de $\Z_K$ d\'efini par
$$
\Delta_{E/K}^{\min}=\prod_{v\in\Pl^\fin_K}\fp_v^{n_v}.
$$
La norme $\Nm\Delta_{E/K}^{\min}$ de cet id\'eal satisfait \`a
$$
\log\Nm\Delta_{E/K}^{\min}=
\sum_{v\in\Pl^\fin_K}n_v\log\#k_v.
$$
On d\'efinit
$$
\Delta_E^\stable=(\Nm\Delta_{E/K}^{\min})^{1/[K:\Q]}.
$$
On note que $\Delta_E^\stable$ ne d\'epend pas du choix de~$K$ et
peut \^etre calcul\'e \`a partir de la factorisation (ou l'id\'eal
d\'enominateur) de~$j_E$ dans n'importe quel corps de nombres
contenant $j_E$; il n'est pas n\'ecessaire de conna{\^\i}tre un
corps~$K$ sur lequel $E$ a r\'eduction semi-stable.

\subsection Le th\'eor\`eme principal

\proclaimx Th\'eor\`eme (cf.\ Cremona, Prickett et
Siksek \citex{Cremona-Prickett-Siksek}{Theorem~1}). Soit $E$ une
courbe elliptique sur~$\Qbar$ donn\'ee par une \'equation de
Weierstra\ss\ \`a coefficients alg\'ebriquement entiers.  Alors pour
tout corps de nombres $K\subset\Qbar$ tel que $E$ soit d\'efinie
sur~$K$, on a
$$
\eqalign{
[K:\Q]\inf_{P\in E(\Qbar)}(h(P)-\hat h(P))&=
\sum_{v\in\Pl^\infin_K}\epsilon_v\inf_{E(\bar K_v)}\phi_v
-{1\over6}\log|\Norm_{K/\Q}\Delta_E|,\cr
[K:\Q]\sup_{P\in E(\Qbar)}(h(P)-\hat h(P))&=
\sum_{v\in\Pl^\infin_K}\epsilon_v\sup_{E(\bar K_v)}\phi_v
+{[K:\Q]\over12}\log\Delta_E^\stable.}
$$

\label{theoreme}

\preuve Quitte \`a \'elargir $K$, on peut supposer que $K$ soit
totalement complexe, que $E$ ait r\'eduction semi-stable scind\'ee
sur~$K$ et que tous les $n_v$ soient pairs.  Pour chaque place finie
$v$ de~$K$, l'hypoth\`ese que le mod\`ele de Weierstra\ss\ soit entier
par rapport \`a~$v$ implique
$$
\eqalign{
\inf_{E(K_v)}\phi_v&={1\over6\epsilon_v}\log|\Delta_E|_v,\cr
\sup_{E(K_v)}\phi_v&=-{1\over12\epsilon_v}\log|\Delta_{E/K}^{\min}|_v;}
\eqnumber{cond-loc}
$$
voir \citex{Cremona-Prickett-Siksek}{Proposition~8}.  (Dans loc.\
cit.\ les r\'esultats sont donn\'es en termes de la fonction
$\Psi_v(P)=\epsilon_v\phi_v(P)+{1\over6}\log|\Delta_E|_v$.)  On en
d\'eduit que
$$
\eqalign{
[K:\Q](h(P)-\hat h(P))&\ge
\sum_{v\in\Pl^\infin_K}\epsilon_v\inf_{E(K_v)}\phi_v
-{1\over6}\log|\Norm_{K/\Q}\Delta_E|,\cr
[K:\Q](h(P)-\hat h(P))&\le
\sum_{v\in\Pl^\infin_K}\epsilon_v\sup_{E(K_v)}\phi_v
+{1\over12}\log\Nm\Delta_{E/K}^{\min}\cr
&=\sum_{v\in\Pl^\infin_K}\epsilon_v\sup_{E(K_v)}\phi_v
+{[K:\Q]\over12}\log\Delta_E^\stable.}
$$

Il reste \`a d\'emontrer que ces bornes {\it localement\/} optimales
donnent des bornes {\it globalement\/} optimales.  \`A cet effet, on
utilise l'approximation sur $\P^1$.  Pour chaque place finie $v$, la
fonction $\phi_v$ atteint son maximum aux points de $E(K_v)$ dont la
$x$-coordonn\'ee est dans un certain ouvert non vide $U_v$ de
$\P^1(K_v)$; de plus, on a $U_v=\P^1(K_v)$ pour toutes les
$v\in\Pl^\fin_K$ sauf un nombre fini.  De fa\c{c}on analogue, pour
chaque place infinie (complexe) $v$ et tout $\epsilon>0$, la fonction
$\phi_v$ sur~$E(K_v)$ est $\epsilon$-proche de son supremum sur
l'ensemble des points dont la $x$-coordonn\'ee est dans un certain
ouvert non vide $U_v$ de $\P^1(K_v)$.  Par approximation, il existe
$x_\epsilon\in\P^1(K)$ dont l'image dans $\P^1(K_v)$ est dans $U_v$
pour toute place~$v$.  Soit $P_\epsilon\in E(\Qbar)$ un point avec
$x(P_\epsilon)=x_\epsilon$.  (On remarquera que $P_\epsilon$ est
d\'efini sur un extension quadratique de~$K$.)  En faisant
$\epsilon\to0$, on obtient une suite de points $P_\epsilon$ pour
lesquels $h(P_\epsilon)-\hat h(P_\epsilon)$ converge vers la borne
sup\'erieure d\'esir\'ee.  Le m\^eme argument marche pour la borne
inf\'erieure.\preuvefinie

Vu le th\'eor\`eme~\ref{theoreme}, il reste \`a \'etudier les
fonctions $\phi_v$ sur $E(\bar K_v)$ pour $v$ une place
archim\'edienne.

\section Pr\'eliminaires sur les r\'eseaux

Soit $\Lambda$ un r\'eseau dans~$\C$.  On note
$$
\vol_\Lambda={i\over2}\int_{\C/\Lambda} dz\wedge\bar dz.
$$
Soit $\mucan_\Lambda$ la $(1,1)$-forme canonique sur~$\C/\Lambda$,
d\'efinie par
$$
\mucan_\Lambda={1\over\vol_\Lambda}{i\over2}dz\wedge d\bar z.
$$

Rappelons la d\'efinition des fonctions $\sigma$, $\zeta$ et $\wp$ de
Weierstra\ss:
$$
\eqalign{
\sigma_\Lambda(z)&=z\prod_{\textstyle{\omega\in\Lambda\atop\omega\ne0}}
\biggl(1-{z\over\omega}\biggr)
\exp\biggl({z\over\omega}+{z^2\over 2\omega^2}\biggr),\cr
\zeta_\Lambda(z)={\sigma_\Lambda'\over\sigma_\Lambda}(z)
&={1\over z}+\sum_{\textstyle{\omega\in\Lambda\atop\omega\ne0}}
\biggl({1\over z-\omega}+{1\over\omega}+{z\over\omega^2}\biggr),\cr
\wp_\Lambda(z)=-\zeta_\Lambda'(z)
&={1\over z^2}+\sum_{\textstyle{\omega\in\Lambda\atop\omega\ne0}}
\biggl({1\over(z-\omega)^2}-{1\over\omega^2}\biggr).}
$$
On rappelle que la fonction $\wp_\Lambda$ est p\'eriodique par rapport
\`a~$\Lambda$, et que la fonction $\zeta_\Lambda$ est
quasi-p\'eriodique: il existe un homomorphisme
$$
\eta_\Lambda:\Lambda\to\C
$$
tel que
$$
\zeta_\Lambda(z+\omega)=\zeta_\Lambda(z)+\eta_\Lambda(\omega)
\quad\hbox{pour tout }\omega\in\Lambda.
$$
Soient $g_2(\Lambda)$ et $g_3(\Lambda)$ les nombres complexes tels que
$$
\wp_\Lambda'(z)^2=
4\wp_\Lambda(z)^3-g_2(\Lambda)\wp_\Lambda(z)-g_3(\Lambda),
$$
et soit 
$$
\Delta_\Lambda=g_2(\Lambda)^3-27g_3(\Lambda)^2.
$$

Soit $Q_\Lambda:\C\to\R$ l'unique forme $\R$-quadratique
satisfaisant \`a
$$
Q_\Lambda(\omega)=\Re(\omega\cdot\eta_\Lambda(\omega))
\quad\hbox{pour tout }\omega\in\Lambda.
$$
On consid\`ere la fonction
$$
\eqalign{
\lambda_\Lambda:\C/\Lambda\setminus\{0\}&\longrightarrow\R\cr
z&\longmapsto-\log|\sigma_\Lambda(z)|+{1\over 2}Q_\Lambda(z)
-{1\over 12}\log|\Delta_\Lambda|.}
$$
Cette fonction est lisse sur $\C/\Lambda\setminus\{0\}$ et poss\`ede
une singularit\'e logarithmique en~$0$.  Elle satisfait \`a
l'\'equation diff\'erentielle
$$
2i\partial\bar\partial\lambda_\Lambda=
2\pi(\mucan_\Lambda-\delta_0)
$$
avec la normalisation
$$
\int_{\C/\Lambda}\lambda_\Lambda\mucan_\Lambda=0.
$$
Comme $\lambda_\Lambda$ prend des valeurs r\'eelles, il existe une
unique fonction lisse (mais non holomorphe)
$Z_\Lambda:\C/\Lambda\setminus\{0\}\to\C$ telle que
$$
d\lambda_\Lambda(z)=-{1\over2}
\bigl(Z_\Lambda(z)dz+\overline{Z_\Lambda(z)}d\bar z\bigr).
\eqnumber{eq-dlambda}
$$

Il existe $C_\Lambda\in\C$ et $D_\Lambda\in\R$ tels que
$$
Q_\Lambda(z)={C_\Lambda\over 2} z^2+{\bar C_\Lambda\over 2}\bar z^2
+D_\Lambda z\bar z
$$
et donc
$$
Z_\Lambda(z)=\zeta_\Lambda(z)-C_\Lambda z-D_\Lambda\bar z.
$$

Soit $(\omega_1,\omega_2)$ une $\Z$-base de~$\Lambda$ avec
$\Im(\omega_1/\omega_2)>0$.  Alors on a
$$
\vol_\Lambda=\Im(\omega_1\bar\omega_2).
$$
On pose
$$
\eqalign{
\eta_1&=\eta_\Lambda(\omega_1)=2\zeta_\Lambda(\omega_1/2),\cr
\eta_2&=\eta_\Lambda(\omega_2)=2\zeta_\Lambda(\omega_2/2).}
$$
En utilisant la relation de Legendre
$$
\eta_2\omega_1-\eta_1\omega_2=2\pi i,
$$
il est facile de montrer que
$$
\eqalign{
C_\Lambda&={\eta_1\bar\omega_2-\eta_2\bar\omega_1\over 2i\vol_\Lambda},\cr
D_\Lambda&={\pi\over\vol_\Lambda}.}
$$

\section \'Etude des fonctions $\phi_v$ aux places archim\'ediennes

Soit $K$ un corps de nombres, et soit $E$ une courbe elliptique
sur~$K$ donn\'ee par une \'equation de
Weierstra\ss~\eqref{weierstrass}.  Soit $v$ une place archim\'edienne
de~$K$.  On fixe un plongement $K\to\C$ correspondant \`a~$v$, et on
regarde $E$ comme courbe elliptique sur~$\C$ au moyen de ce
plongement.

Soit $\Lambda_v\subset\C$ le r\'eseau des p\'eriodes de~$E$ par
rapport \`a la $1$-forme standard $dx/(2y+a_1x+a_3)$ du mod\`ele de
Weierstra\ss\ donn\'e.  On note $\wp_v$ la fonction de Weierstra\ss\
relative au r\'eseau~$\Lambda_v$, vue comme fonction m\'eromorphe
sur~$\C/\Lambda_v$.  Pour $P\in E(\C)$, on note $z_P$ le point
correspondant de $\C/\Lambda_v$.  Alors on a
$$
\wp_v(z_P)=x(P)+{b_2\over12}.
$$
De plus, on a
$$
g_2(\Lambda_v)={c_4\over 12},\quad
g_3(\Lambda_v)={c_6\over 216}.
$$

On note
$$
\lambda_v(z)=\lambda_{\Lambda_v}(z).
$$
La fonction~$\lambda_v$ est la hauteur locale \`a la place~$v$; voir
Silverman \citex{Silverman, ATAEC}{Theorem VI.3.2}.

La fonction $\phi_v$ est donn\'ee par
$$
\phi_v(z)=\log\max\left\{1,\left|\wp_v(z)-{b_2\over12}\right|\right\}
-2\lambda_v(z)
\quad\hbox{pour tout }z\in\C/\Lambda_v.
$$
Soient $t_1,t_2\in\C/\Lambda_v$ les z\'eros de la
fonction~$\wp_v(z)-{b_2\over12}$; on a $t_1+t_2=0$.  Alors
les deux fonctions $\log\left|\wp_v(z)-{b_2\over12}\right|$ et
$2\lambda_v(z)-\lambda_v(z-t_1)-\lambda_v(z-t_2)$ ont m\^eme image
sous l'op\'erateur laplacien, \`a savoir
$2\pi(\delta_{t_1}+\delta_{t_2}-2\delta_0)$.  Il s'ensuit que
$$
\log\biggl|\wp_v(z)-{b_2\over12}\biggr|=
2\lambda_v(z)-\lambda_v(z-t_1)-\lambda_v(z-t_2)+I_v,
$$
o\`u
$$
I_v=\int_{\C/\Lambda_v}
\log\biggl|\wp_v(z)-{b_2\over12}\biggr|\mucan_{\Lambda_v}.
$$
Cela implique
$$
\phi_v(z)=\cases{
-\lambda_v(z-t_1)-\lambda_v(z-t_2)+I_v&
  si $\left|\wp_v(z)-{b_2\over12}\right|\ge1$,\cr
-2\lambda_v(z)&
  si $\left|\wp_v(z)-{b_2\over12}\right|\le1$.}
\eqnumber{eq-phi}
$$
On pose
$$
S=\left\{z\in\C/\Lambda_v\biggm|
\left|\wp_v(z)-{b_2\over12}\right|=1\right\}.
$$
On peut utiliser \eqref{eq-phi} pour calculer $\phi_v(z)$.  La
constante $I_v$ peut \^etre d\'etermin\'ee en comparant les deux
expressions pour $\phi_v(z)$ pour n'importe quel $z\in S$.

En utilisant les formules \eqref{eq-phi} et~\eqref{eq-dlambda}, on
voit que la d\'eriv\'ee de $\phi_v$ est donn\'ee par
$$
d\phi_v(z)=W_v(z)dz+\overline{W_v(z)}d\bar z
\quad\hbox{pour }\left|\wp_v(z)-{b_2\over12}\right|\ne 1,
\eqnumber{dphi}
$$
o\`u $W_v$ est la fonction continue sur $\C/\Lambda_v\setminus S$
d\'efinie par
$$
W_v(z)=\cases{
{1\over2}\bigl(Z_{\Lambda_v}(z-t_1)+Z_{\Lambda_v}(z-t_2)\bigr)&
  si $\left|\wp_v(z)-{b_2\over12}\right|>1$,\cr
Z_{\Lambda_v}(z)&
  si $\left|\wp_v(z)-{b_2\over12}\right|<1$.}
$$

\proclaim Lemme. La d\'eriv\'ee de $Z_{\Lambda_v}$ est
$$
dZ_{\Lambda_v}(z)=
\bigl(-\wp_v(z)-C_{\Lambda_v}\bigr)dz-D_{\Lambda_v} d\bar z.
$$
La d\'eriv\'ee de
${1\over2}\bigl(Z_{\Lambda_v}(z-t_1)+Z_{\Lambda_v}(z-t_2)\bigr)$ est
$$
d\biggl({1\over2}\bigl(Z_{\Lambda_v}(z-t_1)+Z_{\Lambda_v}(z-t_2)\bigr)\biggr)
=\biggl(-C_{\Lambda_v}-{b_2\over12}-{1\over2}{b_4\over \wp_v(z)-{b_2\over12}}
-{1\over2}{b_6\over(\wp_v(z)-{b_2\over12})^2}\biggr)dz
-D_\Lambda d\bar z.
$$

\label{deriv-W}

\preuve La premi\`ere partie suit de la d\'efinition de $Z_{\Lambda_v}$
et du fait que $\zeta_{\Lambda_v}'(z)=-\wp_v(z)$.

En rempla\c{c}ant $z$ par $z-t$, on obtient
$$
dZ_{\Lambda_v}(z-t)
=\bigl(-\wp_v(z-t)-C_{\Lambda_v}\bigr)dz-D_{\Lambda_v} d\bar z,
$$
de sorte que
$$
{1\over2}d\bigl(Z_{\Lambda_v}(z-t_1)+Z_{\Lambda_v}(z-t_2)\bigr)
=\left(-{\wp_v(z-t_1)+\wp_v(z-t_2)\over2}
-C_{\Lambda_v}\right)dz-D_{\Lambda_v} d\bar z.
$$
La loi d'addition permet d'exprimer $\wp_v(z-t_1)+\wp_v(z-t_2)$ en
$\wp_v(z)$ et $\wp_v(t_1)=\wp_v(t_2)={b_2\over12}$ comme suit:
$$
\eqalign{
\wp_v(z-t_1)+\wp_v(z-t_2)&=
{4{b_2\over12}\wp_v(z)(\wp_v(z)+{b_2\over12})
-g_2(\Lambda_v)(\wp_v(z)+{b_2\over12})
-2g_3(\Lambda_v)\over 2(\wp_v(z)-{b_2\over12})^2}\cr
&={b_2\over6}+{12({b_2\over12})^2-g_2(\Lambda_v)
\over 2(\wp_v(z)-{b_2\over12})}
+{4({b_2\over12})^3-g_2(\Lambda_v){b_2\over12}-g_3(\Lambda_v)
\over(\wp_v(z)-{b_2\over12})^2}\cr
&={b_2\over6}+{12({b_2\over12})^2-{c_4\over12}
\over 2(\wp_v(z)-{b_2\over12})}
+{4({b_2\over12})^3-{c_4\over12}{b_2\over12}-{c_6\over216}
\over (\wp_v(z)-{b_2\over12})^2}\cr
&={b_2\over6}+{b_4\over \wp_v(z)-{b_2\over12}}
+{b_6\over(\wp_v(z)-{b_2\over12})^2}.}
$$
Ceci implique la deuxi\`eme partie.\preuvefinie

On note que $W_v$ est harmonique, de sorte que pour borner $W_v$ sur
$\C/\Lambda_v\setminus S$, il suffit de borner les fonctions
$Z_{\Lambda_v}(z)$ et
${1\over2}\bigl(Z_{\Lambda_v}(z-t_1)+Z_{\Lambda_v}(z-t_2)\bigr)$
sur~$S$.  Pour chacune de ces deux fonctions, le lemme suivant donne
une majoration de sa valeur absolue sur~$S$ \`a partir d'une
seule \'evaluation et d'une majoration de quelques nombres r\'eels
associ\'es \`a~$E$.

\def\e#1{\exp(#1)}

\proclaim Lemme. Pour tout $p,q\in S$, on a
$$
\eqalign{
\bigl|Z_{\Lambda_v}(q)\bigr|&\le\bigl|Z_{\Lambda_v}(p)\bigr|+M_1J,\cr
\left|{1\over2}\bigl(Z_{\Lambda_v}(q-t_1)+Z_{\Lambda_v}(q-t_2)\bigr)\right|&\le
\left|{1\over2}\bigl(Z_{\Lambda_v}(p-t_1)+Z_{\Lambda_v}(p-t_2)\bigr)\right|
+M_2J,}
$$
o\`u
$$
\eqalign{
M_1&=\left|C_{\Lambda_v}+{b_2\over12}\right|+|D_{\Lambda_v}|+1,\cr
M_2&=\left|C_{\Lambda_v}+{b_2\over12}\right|+|D_{\Lambda_v}|
+{|b_4|\over2}+{|b_6|\over2},\cr
J&=\int_0^{2\pi}{d\theta\over|4\e{3i\theta}+b_2\e{2i\theta}
+2b_4\e{i\theta}+b_6|^{1/2}}.}
$$

\label{lemme-W1}

\preuve Quitte \`a remplacer $p$ par $-p$ (ce qui a pour effet de
multiplier $Z_{\Lambda_v}(p)$ par $-1$), on peut supposer qu'il existe
un chemin $\gamma$ de $p$ \`a $q$ dans~$S$ tel que la fonction
$\wp_v(z)-{b_2\over12}$ identifie $\gamma$ avec un segment du cercle
unit\'e dans~$\C$.  On a
$$
\eqalign{
Z_{\Lambda_v}(q)&=Z_{\Lambda_v}(p)+\int_\gamma d Z_{\Lambda_v}\cr
&=Z_{\Lambda_v}(p)-\int_\gamma \bigl((\wp_v(z)+C_{\Lambda_v})dz
+D_{\Lambda_v} d\bar z\bigr).}
$$
Pour tout $z\in S$, le lemme \ref{deriv-W} implique
$$
\eqalign{
|\wp_v(z)+C_{\Lambda_v}|&\le\left|\wp_v(z)-{b_2\over12}\right|
+\left|C_{\Lambda_v}+{b_2\over12}\right|\cr
&=1+\left|C_{\Lambda_v}+{b_2\over12}\right|.}
$$
On en d\'eduit que
$$
|Z_{\Lambda_v}(q)|\le|Z_{\Lambda_v}(p)|+M_1\int_\gamma |dz|.
$$
De fa\c{c}on analogue, on obtient
$$
\left|{1\over2}\bigl(Z_{\Lambda_v}(q-t_1)+Z_{\Lambda_v}(q-t_2)\bigr)\right|\le
\left|{1\over2}\bigl(Z_{\Lambda_v}(p-t_1)+Z_{\Lambda_v}(p-t_2)\bigr)\right|
+M_2\int_\gamma |dz|.
$$
En utilisant la formule
$$
(2y+a_1x+a_3)^2=4x^3 + b_2 x^2 + 2b_4 x + b_6,
$$
on montre facilement que
$$
\eqalign{
\int_\gamma|dz|&\le\int_{|x|=1}{|dx|\over|2y+a_1x+a_3|}\cr
&=\int_0^{2\pi}{d\theta\over|4\e{3i\theta}+b_2\e{2i\theta}
+2b_4\e{i\theta}+b_6|^{1/2}}\cr
&=J.}
$$
Ceci compl\`ete la d\'emonstration.\preuvefinie

\proclaim Corollaire. Soit $p$ un point quelconque de~$S$.
Pour tout $z\in\C/\Lambda_v\setminus S$, on a
$$
|W_v(z)|\le\max\left\{\bigl|Z_{\Lambda_v}(p)\bigr|+M_1J,
\left|{1\over2}\bigl(Z_{\Lambda_v}(p-t_1)+Z_{\Lambda_v}(p-t_2)\bigr)\right|
+M_2J\right\}.
$$

\label{cor-W1}

Les r\'esultats suivants ne sont pas strictement n\'ecessaires, mais
permettent d'obtenir un algorithme plus efficace ci-dessous.

\proclaim Lemme. Soit $R$ un sous-ensemble convexe de~$\C$.
\smallskip
\item{(a)} Si $|\wp_v(z)-{b_2\over12}|<1$ pour tout $z\in R$, on a
$$
|W_v(z)-W_v(z_0)|\le
\left(\left|C_{\Lambda_v}+{b_2\over12}\right|+|D_{\Lambda_v}|+1\right)|z-z_0|
\quad\hbox{pour tout }z,z_0\in R.
$$
\item{(b)} Si $|\wp_v(z)-{b_2\over12}|>1$ pour tout $z\in R$, on a
$$
|W_v(z)-W_v(z_0)|\le
\left(\left|C_{\Lambda_v}+{b_2\over12}\right|+|D_{\Lambda_v}|+
{|b_4|+|b_6|\over2}\right)|z-z_0|
\quad\hbox{pour tout }z,z_0\in R.
$$

\label{lemme-W2}

\preuve Dans chacun des deux cas, la fonction $W_v$ est
diff\'erentiable sur~$R$, et il suffit de borner sa d\'eriv\'ee.  Le
lemme \ref{deriv-W} implique que pour $|\wp_v(z)-{b_2\over 12}|\le 1$,
on a
$$
\eqalign{
\left|\wp_v(z)+C_{\Lambda_v}\right|
&\le\left|\wp_v(z)-{b_2\over12}\right|
+\left|C_{\Lambda_v}+{b_2\over12}\right|\cr
&\le 1+\left|C_{\Lambda_v}+{b_2\over12}\right|.}
$$
Pour $|\wp_v(z)-{b_2\over12}|\ge 1$, le lemme \ref{deriv-W} implique
$$
\left|{\wp_v(z-t_1)+\wp_v(z-t_2)\over2}+C_{\Lambda_v}\right|
\le\left|C_{\Lambda_v}+{b_2\over12}\right|
+{|b_4|\over2}+{|b_6|\over2},
$$
ce qui compl\`ete la d\'emonstration.\preuvefinie

On consid\`ere maintenant des parall\'elogrammes de la forme
$$
R(z_0,z_1,z_2)=\left\{z_0+s_1z_1+s_2z_2\bigm|
s_1,s_2\in\left[-\textfrac1/2,\textfrac1/2\right]\right\}
$$
pour $z_0,z_1,z_2\in\C$ tels que $z_1$ et~$z_2$ sont
$\R$-lin\'eairement ind\'ependants.  
On note
$$
\eqalign{
d(z_1,z_2)&=\sup_{z\in R(z_0,z_1,z_2)}|z-z_0|\cr
&={1\over2}\max\{|z_1-z_2|,|z_1+z_2|\}.}
$$

\proclaim Corollaire. Soit $R=R(z_0,z_1,z_2)$ comme ci-dessus.
\smallskip
\item{(a)} Si $|\wp_v(z)-{b_2\over12}|<1$ pour tout $z\in R$, on a
$$
|W_v(z)|\le|W_v(z_0)|+
\left(\left|C_{\Lambda_v}+{b_2\over12}\right|+|D_{\Lambda_v}|+1\right)
d(z_1,z_2)
\quad\hbox{pour tout }z\in R.
$$
\item{(b)} Si $|\wp_v(z)-{b_2\over12}|>1$ pour tout $z\in R$, on a
$$
|W_v(z)|\le|W_v(z_0)|+
\left(\left|C_{\Lambda_v}+{b_2\over12}\right|+|D_{\Lambda_v}|+
{|b_4|+|b_6|\over2}\right)d(z_1,z_2)
\quad\hbox{pour tout }z\in R.
$$

\label{cor-W2}

\section Un algorithme

L'\'etude de $\phi_v$ ci-dessus permet de construire un algorithme
pour calculer le supremum de $\phi_v$ avec une pr\'ecision pr\'escrite
$\epsilon$.  On a un algorithme compl\`etement analogue pour calculer
l'infimum.

Notre algorithme, dont l'id\'ee fondamentale est inspir\'ee de
l'algorithme de Cremona, Prickett et
Siksek \citex{Cremona-Prickett-Siksek}{\S9}, fonctionne par dichotomie
r\'ecursive.  On commence avec un domaine fondamental
$R(0,\omega_1,\omega_2)$, o\`u $(\omega_1,\omega_2)$ est une $\Z$-base
de $\Lambda_v$, et on pose $\mu=\phi_v(0)$.  \`A chaque \'etape, on
consid\`ere un parall\'elogramme $R(z_0,z_1,z_2)$.  On remplace $\mu$
par $\max\{\mu,\phi_v(z_0)\}$, de sorte que $\mu$ est toujours la plus
grande valeur de~$\phi_v$ qu'on a rencontr\'e jusque-l\`a.  De plus,
on calcule un majorant $M$ pour la fonction $|W_v|$
sur~$R(z_0,z_1,z_2)$ par le corollaire~\ref{cor-W1} ou le
corollaire~\ref{cor-W2}.  Pour tout $z\in R(z_0,z_1,z_2)$, on a
par~\eqref{dphi}
$$
\eqalign{
|\phi_v(z)-\phi_v(z_0)|&\le
2d(z_1,z_2)\sup_{R(z_0,z_1,z_2)}|W_v|\cr
&\le 2d(z_1,z_2)M.}
\eqnumber{iter}
$$
Dans le cas o\`u
$$
\phi_v(z_0)+2d(z_1,z_2)M<\mu+\epsilon,
$$
on conclut gr\^ace \`a~\eqref{iter} que $\phi_v$ est
inf\'erieur \`a~$\mu+\epsilon$ sur tout le parall\'elogramme
$R(z_0,z_1,z_2)$.  Dans le cas oppos\'e, on coupe $R(z_0,z_1,z_2)$ en
deux nouveaux parall\'elogrammes le long de la droite qui passe par le
centre $z_0$ et qui est parall\`ele \`a un c\^ot\'e de longueur
minimale de~$R(z_0,z_1,z_2)$, et on applique le processus de fa\c{c}on
r\'ecursive \`a ces nouveaux parall\'elogrammes.  Il est facile de
voir que cet algorithme termine et produit un $\mu$ qui satisfait \`a
$$
\sup_{\C/\Lambda_v}\phi_v-\epsilon<\mu\le\sup_{\C/\Lambda_v}\phi_v.
$$

\section Exemples

Voici quelques exemples.  Nous avons utilis\'e {\sc pari/gp}
\cite{pari} pour la plupart des calculs, et Sage \cite{sage} pour
calculer la hauteur canonique de points d\'efinis sur des corps de
nombres.

\subsection La courbe {\rm 11A3}

Cette courbe a un mod\`ele globalement minimal donn\'e par
l'\'equation
$$
E : y^2 + y = x^3 - x^2.
$$
On a
$$
-\Delta_E=\Delta_E^\stable=11.
$$
Notre algorithme donne
$$
\eqalign{
\sup_{E(\C)}\phi_v&=0.597\ldots,\cr
\inf_{E(\C)}\phi_v&=-0.156\ldots}
$$
En utilisant le th\'eor\`eme \ref{theoreme}, on obtient
$$
-0.556 < h(P)-\hat h(P) < 0.798
\quad\hbox{pour tout }P\in E(\Qbar).
$$
Pour comparaison, la majoration trouv\'ee par l'algorithme de
Silverman~\cite{Silverman} est
$$
h(P)-\hat h(P)<4.695,
$$
et celle trouv\'ee par l'algorithme de Cremona, Prickett et
Siksek \cite{Cremona-Prickett-Siksek} pour les $\Q$-points est
$$
h(P)-\hat h(P)<0.300.
$$

Par approximation dans $\P^1$ comme dans la preuve du
th\'eor\`eme~\ref{theoreme}, on peut trouver des points $P,Q\in
E(\Qbar)$ tels que $h(P)-\hat h(P)$ est tr\`es proche de $-0.556$ et
$h(Q)-\hat h(Q)$ est tr\`es proche de $0.798$.  On prend
$$
\eqalign{
P&=(-1,\alpha)\quad\hbox{avec }\alpha^2+\alpha+2=0,\cr
Q&=(37/61,\beta)\quad\hbox{avec }\beta^2+\beta
+{2^3\cdot 3\cdot 37^2\over 61^3}=0.}
$$
Le point $P$ est d\'efini sur $\Q(\sqrt{-7})$; le point $Q$ est
d\'efini sur $\Q(\sqrt{7\cdot 11\cdot 17\cdot 61\cdot 73})$.  On a
$$
\vcenter{\halign{
$#$\hfil\quad&$#$\hfil\quad&$#$\hfil\cr
h(P)=0,&
\hat h(P)=0.5556807\ldots,&
h(P)-\hat h(P)=-0.5556807\ldots,\cr
h(Q)=\log 61,&
\hat h(Q)=3.3130740\ldots,&
h(Q)-\hat h(Q)=0.7977997\ldots\cr}}
$$

\subsection La courbe {\rm 15A4}

Voici un example pour montrer que notre m\'ethode donne parfois une
meilleure majoration de $h-\hat h$ pour les $\Qbar$-points que celle
de Cremona, Prickett et Siksek \cite{Cremona-Prickett-Siksek} pour les
$\Q$-points.

La courbe 15A4 a un mod\`ele globalement minimal donn\'e par
l'\'equation
$$
E : y^2 + xy + y = x^3 + x^2 + 35x - 28.
$$
On a
$$
-\Delta_E=\Delta_E^\stable=3^2\cdot 5^8.
$$
Notre algorithme donne
$$
\eqalign{
\inf_{E(\C)}\phi_v&=0.584\ldots,\cr
\sup_{E(\C)}\phi_v&=2.512\ldots}
$$
En utilisant le th\'eor\`eme \ref{theoreme}, on obtient
$$
-1.928 < h(P)-\hat h(P) < 3.769
\quad\hbox{pour tout }P\in E(\Qbar).
$$
Pour comparaison, la majoration trouv\'ee par l'algorithme
de \cite{Cremona-Prickett-Siksek} est
$$
h(P)-\hat h(P)<3.915
\quad\hbox{pour tout }P\in E(\Q).
$$

\subsection La courbe {\rm 5077A1}

Cette courbe, d\'ej\`a \'etudi\'ee par Buhler, Gross et Zagier
\cite{Buhler-Gross-Zagier}, a un mod\`ele globalement minimal donn\'e
par l'\'equation
$$
E : y^2 + y = x^3 - 7x + 6
$$
On a
$$
\Delta_E=\Delta_E^\stable=5077.
$$
Notre algorithme donne
$$
\eqalign{
\inf_{E(\C)}\phi_v&=0.217\ldots,\cr
\sup_{E(\C)}\phi_v&=1.422\ldots}
$$
En utilisant le th\'eor\`eme \ref{theoreme}, on obtient
$$
-1.206 < h(P)-\hat h(P) < 2.134
\quad\hbox{pour tout }P\in E(\Qbar).
$$
Des bornes optimales pour les $\Q$-points ont d\'ej\`a \'et\'e
calcul\'ees dans loc.\ cit\null.  En fait, on a
$$
-1.2050811\ldots\le h(P)-\hat h(P)\le 0
\quad\hbox{pour tout }P\in E(\Q),
$$
o\`u la borne inf\'erieure est atteinte par le point $(-1,3)$ et la
borne sup\'erieure par le point \`a l'infini.

Par approximation dans $\P^1$, on trouve le point
$$
P=(5169,\alpha)\quad\hbox{avec }\alpha^2+\alpha
-138108205632 
=0,
$$
d\'efini sur $\Q(\sqrt{7\cdot 5077\cdot 15544411})$.  On a
$$
h(P)=\log 5169,\quad
\hat h(P)=6.4174217\ldots,\quad
h(P)-\hat h(P)=2.1330128\ldots
$$

\subsection Une courbe \'etudi\'ee par Cremona, Prickett et Siksek

Notre dernier exemple est une courbe elliptique de rang~4
sur~$\Q$ \'etudi\'ee par Cremona, Prickett et Siksek
dans \citex{Cremona-Prickett-Siksek}{\S11}:
$$
E : y^2 =  x^3 - 459 x^2 - 3478 x + 169057.
$$
Cette courbe n'est pas semi-stable sur~$\Q$, le conducteur \'etant
$2^2\cdot 199\cdot 362793983647$.  On a
$$
\eqalign{
\Delta_E&=2^4\cdot 199\cdot 362793983647,\cr
\Delta_E^\stable&=199\cdot 362793983647.}
$$
Notre algorithme donne
$$
\eqalign{
\inf_{E(\C)}\phi_v&=0.879\ldots,\cr
\sup_{E(\C)}\phi_v&=5.780\ldots}
$$
En utilisant le th\'eor\`eme \ref{theoreme}, on obtient
$$
-4.901 < h(P)-\hat h(P) < 8.440
\quad\hbox{pour tout }P\in E(\Qbar).
$$
Les bornes trouv\'ees par Cremona, Prickett et
Siksek \cite{Cremona-Prickett-Siksek} sont
$$
-6.532 < h(P)-\hat h(P) < 0.4621
\quad\hbox{pour tout }P\in E(\Q).
$$
Ils ont \'egalement trouv\'e un point $P\in E(\Q)$ pour lequel
$$
h(P)-\hat h(P)=-4.9001533\ldots
$$
De l'autre c\^ot\'e, le point
$$
Q=(-45092013952912,\alpha)
\quad\hbox{avec }
\alpha^2=-199\cdot 1601\cdot 22133
\cdot 362793983647\cdot 35838855272124651419
$$
satisfait \`a
$$
h(Q)=31.4397262\ldots,\quad
\hat h(Q)=23.0000267\ldots,\quad
h(Q)-\hat h(Q)=8.4396995\ldots
$$

\references{Bibliographie}

\NoAutoSpaceBeforeFDP

\def\No{{\cmcyr\char"19}}

\reference{Buhler-Gross-Zagier} J. P. {\sc Buhler}, B. H. {\sc Gross}
and D. B. {\sc Zagier}, On the conjecture of Birch and Swinnerton-Dyer
for an elliptic curve of rank~3.  {\it Mathematics of Computation}
{\bf 44} (1985), 473--481.

\reference{Cremona-Prickett-Siksek} J. E. {\sc Cremona}, M. {\sc
  Prickett} and S. {\sc Siksek}, Height difference bounds for elliptic
curves over number fields.  {\it Journal of Number Theory\/} {\bf 116}
(2006), no.~1, 42--68.

\reference{Demjanenko} {\cyr\V}.~{\cyr\A}.\
{\cycsc\D\ye\m\soft\ya\n\ye\n\k\o}, {\cyr\O\ts\ye\n\k\a\
  \o\s\t\a\t\o\ch\n\o\g\o\ \ch\l\ye\n\a\ \v\ \f\o\r\m\u\l\ye\
  \T\e\j\t\a}. {\cyti\M\a\t\ye\m\a\t\i\ch\ye\s\k\i\ye\
  \Z\a\m\ye\t\k\i\/} {\bf 3} (1968), \No~3, 271--278.\hfill\break
V. A. {\sc Dem$'$janenko}, An estimate of the
remainder term in Tate's formula.  {\it Mathematical Notes\/} {\bf 3}
(1968), no.~3, 173--177.  (English translation.)

\reference{pari} The PARI Group, {\sc pari/gp}, version 2.6.0.
Bordeaux, 2012,\hfill\break {\tt http://pari.math.u-bordeaux.fr/}.

\reference{Siksek} S. {\sc Siksek}, Infinite descent on elliptic
curves.  {\it Rocky Mountain Journal of Mathematics\/} {\bf 25}
(1995), no.~4, 1501--1538.

\reference{Silverman} J. H. {\sc Silverman}, The difference between
the Weil height and the canonical height on elliptic curves.  {\it
  Mathematics of Computation\/} {\bf 55} (1990), 723--743.

\reference{Silverman, ATAEC} J. H. {\sc Silverman}, {\sl Advanced
  Topics in the Arithmetic of Elliptic Curves\/}.  Graduate Texts in
Mathematics {\bf 151}.  Springer-Verlag, New York, 1994.

\reference{sage} W. A. {\sc Stein} et~al., Sage Mathematics Software,
version 5.2.  The Sage Development Team, 2012, {\tt
  http://www.sagemath.org/}.

\reference{Tate} J. T. {\sc Tate}, The arithmetic of elliptic curves.
{\it Inventiones mathematicae\/} {\bf 23} (1974), 179--206.

\reference{Uchida} Y. {\sc Uchida}, The difference between the
ordinary height and the canonical height on elliptic curves.  {\it
  Journal of Number Theory\/} {\bf 128} (2008), no.~2, 263--279.

\reference{Zimmer} H. G. {\sc Zimmer}, On the difference between the
Weil height and the N\'eron--Tate height.  {\it Mathe\-ma\-ti\-sche
  Zeitschrift\/} {\bf 147} (1976), 35--51.

\vskip2cm
\vbox{
\leftline{Peter Bruin}
\leftline{Institut f\"ur Mathematik}
\leftline{Universit\"at Z\"urich}
\leftline{Winterthurerstrasse 190}
\leftline{CH-8057 Z\"urich}
\smallskip
\leftline{\tt peter.bruin@math.uzh.ch}}

\bye